\documentclass[12pt]{amsart}
\usepackage{a4wide,enumerate,color,graphicx}
\usepackage{amsmath}

 \newtheorem{definition}{Definition}
 \newtheorem{thm}{Theorem}
 
 \newtheorem{lem}[thm]{Lemma}
 
 \theoremstyle{definition}
 
 \theoremstyle{remark}

 \DeclareMathOperator{\dive}{div}

\begin{document}
\title[]{Global existence and uniqueness analysis of cross-diffusion multispecies chemotaxis system with volume-filling}

\author[]{Xi Lin}
\address{{Department of Mathematics and Physics, Guangzhou Maritime University,
Guangzhou 510725, Guangdong Province, China}}
\email{linxi@gzmtu.edu.cn}

\address{}
\email{}

\subjclass{Primary 35K40, 35K55; Secondary 35Q92}
\keywords{multispecies chemotaxis system, entropy method, distance functional, cross-diffusion system with volume-filling}





\begin{abstract}
The system of multispecies chemotaxis equations is a cross-diffusion system with volume-filling. In this work, we show that a weak solution of the two species chemotaxis system exists. The entropy method is a major tool in existence analysis of cross-diffusion systems. Previous investigations indicate that traditional entropies in volume-filling cases may not be able to provide required gradient estimates. In this situation, we upgrade existing matrix computation methods to derive gradient estimates. Due to the cross-diffusion phenomenon, the uniqueness of the weak solution to a cross-diffusion system is very difficult to prove in general. In this work, we apply the distance functional to show that when parameters of the chemotaxis system are identical, the weak solution is unique.
\end{abstract}
\maketitle

\section{Introduction}

\subsection{Motivation} In this paper, we analyze a two species cross-diffusion system with volume-filling. The system is called the two species chemotaxis equations, which is used to describe the movement of two different cell types in crowded and uncrowded tissues. Please refer to \cite{1} for a detailed explanation of this model.

In \cite{2}, the two species chemotaxis system was investigated as a cross-diffusion system with volume-filling, but it did not show that a weak solution of this system exists.

Global existence analysis of this chemotaxis system in \cite{2} indicated that how to derive efficient gradient estimates is the key factor for the existence proof. This work is a continuation of these references. We show that a weak solution of the two species chemotaxis system exists.

Let us denote $\mathcal{D}=\{u=(u_1,u_2):u_1,u_2>0,u_1+u_2<1\}$, and $\overline{\mathcal{D}}$ refers to the closure set of $\mathcal{D}$. For $u=(u_1,u_2):\Omega\times(0,\infty)\to\mathcal{D}$, the two species chemotaxis system refers to parabolic partial differential equations in the form of
\begin{equation}\label{equation}
\partial_tu=\dive\big(A(u)\nabla u\big)\text{ in }Q_T=\Omega\times(0,T),\text{ }T>0,
\end{equation}
where $\Omega\subset\mathbb{R}^{d}$, $d\geq 1$ is a bounded domain, and $u_1,u_2$ refer to cell densities. The diffusion matrix $A(u)=(A_{ij}(u))$ is given by
\begin{equation}\label{diffusiontwo}
\begin{aligned}
A(u)&=\left(\begin{array}{cccccc}
\alpha_1(1-u_2)+u_2 & (\alpha_1-1)u_1 \\
 (\alpha_2-1)u_2& \alpha_2(1-u_1)+u_1 \\
\end{array}\right),\quad\alpha_1,\alpha_2>0,\text{ }\alpha_1,\alpha_2\neq 1.
\end{aligned}
\end{equation}

Equations subject to initial and no-flux boundary conditions
\begin{equation}\label{initial}
\begin{aligned}
&(A_{11}(u)\nabla u_1+A_{12}(u)\nabla u_2)\cdot\nu=0,\quad(A_{21}(u)\nabla u_1+A_{22}(u)\nabla u_2)\cdot\nu=0 \text{ }\text{   on  } \partial{\Omega},\text{ } t>0,\\&u_1(0)=u^0_1,u_2(0)=u^0_2 \text{  in  } \Omega,
\end{aligned}
\end{equation}
where $\nu$ is an exterior unit normal vector to $\partial\Omega$ and $u^0=(u^0_1,u^0_2)$ is the initial datum. We show that once $u^0\in\mathcal{D}$, then a weak solution $u\in\overline{\mathcal{D}}$ of (\ref{equation})-(\ref{initial}) exists. 

If the weak solution $u=(u_1,...,u_n)$ of a cross-diffusion system satisfies that $u\in\overline{\mathcal{V}}$, $\mathcal{V}=\{u=(u_1,...,u_n):u_i>0,\sum_{i=1}^{n}u_i<1\}$, then we call the system has the volume-filling effect. Chemotaxis system is a typical cross-diffusion system with volume-filling.

Cross-diffusion systems were divided into two groups in \cite{2}. Degenerate systems \cite{3,4} and Maxwell-Stefan equations \cite{5,6} are typical cross-diffusion systems with volume-filling. The Shigesada-Kawasaki-Teramoto type population system \cite{7,8} and semiconductor model with electron-hole scattering \cite{9,10} are cross-diffusion systems without volume-filling.

After these preparations, we present the definition of weak solutions to (\ref{equation})-(\ref{initial}) and Theorem \ref{mainresult}.
\begin{definition}\label{defweak}
For every $T>0$, $u^0=(u^0_1,u^0_2)\in\mathcal{D}$, if $u=(u_1,u_2)$ satisfies that $u\in\overline{\mathcal{D}}$, and for every test function $\phi=(\phi_1,\phi_2)$, $\phi_1,\phi_2\in L^2(0,T;H^1(\Omega))$,
\begin{equation}\label{chemdef}
\begin{aligned}
&\int_{0}^{T}\langle\partial_tu_1,\phi_1\rangle dt+\int_{Q_T}(\alpha_1\nabla u_1+(\alpha_1-1)(u_1\nabla u_2-u_2\nabla u_1))\cdot\nabla\phi_1dxdt=0,\\&\int_{0}^{T}\langle\partial_tu_2,\phi_2\rangle dt+\int_{Q_T}(\alpha_2\nabla u_2-(\alpha_2-1)(u_1\nabla u_2-u_2\nabla u_1))\cdot\nabla\phi_2dxdt=0,
\end{aligned}
\end{equation}
with the property that $u_i\in L^2(0,T;H^1(\Omega))$, $\partial_tu_i\in L^2(0,T;H^1(\Omega)^{\prime})$, $i=1,2$, then we call $u$ is a weak solution of $(\ref{equation})$-$(\ref{initial})$.
\end{definition}
\begin{thm}\label{mainresult}
Let $u^0=(u^0_1,u^0_2)$ be such that $u^0\in\mathcal{D}$. There exists a weak solution $u\in\overline{\mathcal{D}}$ of $(\ref{equation})$-$(\ref{initial})$.
\end{thm}

How to show the weak solution of a cross-diffusion system is unique is very difficult, due to the existence of cross-diffusion phenomenon. So far few uniqueness results have been derived.

Among these existing uniqueness results, we usually require that the diffusion coefficients, parameters of a cross-diffusion system satisfy certain assumptions. In this work, we will show Theorem \ref{mainresultunique}.
\begin{thm}\label{mainresultunique}
Let $u^0=(u^0_1,u^0_2)$ be such that $u^0\in\mathcal{D}$. The weak solution $u\in\overline{\mathcal{D}}$ of $(\ref{equation})$-$(\ref{initial})$ is unique, if $\alpha_1=\alpha_2=\alpha$.
\end{thm}

Once $\alpha_1\neq\alpha_2$, we conjecture that the weak solution still exists uniquely, though
it remains very difficult to prove this conjecture.

In the last section, we have made some further remarks. We list Maxwell-Stefan equations, degenerate systems and semiconductor model as examples. We try to derive structural relationships among these cross-diffusion systems. We have also indicated the difficulties for the uniqueness proof of cross-diffusion systems.

\subsection{Major Difficulties and Novelties} The entropy method has been demonstrated as a powerful tool in global existence analysis of cross-diffusion systems. How to apply the entropy method to show weak solutions exist to cross-diffusion systems was discussed in \cite{3,4,5,6,7}. Let us present these assumptions:

(H1). There exists a convex function $h\in C^2(\mathcal{V},(0,\infty))$ such that its derivative $h^{\prime}:\mathcal{V}\to\mathbb{R}^n$ is invertible on $\mathbb{R}^n$.

(H2). For every $z=(z_1,...,z_n)^{\mathrm{T}}\in\mathbb{R}^n$ and $u=(u_1,...,u_n)\in\mathcal{V}$, there exists a constant $c>0$ independent of $z$ and $u$, such that
\begin{equation}\nonumber
z^{\mathrm{T}}h^{\prime\prime}(u)A(u)z\geq c\sum_{i=1}^{n}z^2_i.
\end{equation}

In here, $H(u)=(H_{ij}(u))=h^{\prime\prime}(u)$ is the Hessian matrix of $h(u)$, and $A(u)$ is the diffusion matrix of the cross-diffusion system. For every $1\leq i,j\leq n$, $H_{ij}(u)=\partial^2h/\partial u_i\partial u_j$. 

(H3). There exists a constant $c^{\prime}>0$ such that for $u\in\mathcal{V}$ and $i,j=1,...,n$, $|A_{ij}(u)|<c^{\prime}$.

If $h(u)$ and $A(u)$ satisfy assumptions (H1)-(H3), then the entropy method can be applied to show a weak solution of the cross-diffusion system exists. Assumptions (H1)-(H3) are intimately related to the concept of ``entropy structure" in \cite{11}. The entropy structure is an essential structural character of cross-diffusion systems.

In global existence analysis of Maxwell-Stefan equations \cite{5}, the entropy was chosen as
\begin{equation}\label{maxwellentropy}
\begin{aligned}
h(\rho^{\prime})=k\sum_{i=1}^{N+1}x_i(\ln x_i-1)+k,\quad\rho^{\prime}=(\rho_1,...,\rho_N).
\end{aligned}
\end{equation}
 
In (\ref{maxwellentropy}), $x_i$ and $\rho_i$ are related by $x_i=\rho_i/(kM_i)$, $\rho_i>0$, with $k=\sum_{i=1}^{N+1}(\rho_i/M_i)$, and $\sum_{i=1}^{N+1}\rho_i=\sum_{i=1}^{N+1}x_i=1$. $M_i>0$, $i=1,...,N+1$ are constants.

For degenerate systems in \cite{4}, the entropy was
\begin{equation}\label{degenerateentropy}
\begin{aligned}
h(u)=\sum_{i=1}^{n}(u_i\log u_i-u_i+1)+\int_{a}^{u_{n+1}}\log q(s)ds+\chi(u),\quad a\in(0,1],\text{ }\sum_{i=1}^{n+1}u_i=1,\text{ }u_i>0,
\end{aligned}
\end{equation}
and $q(s),\chi(u)$ are functions defined in (7)-(8) of \cite{4}.

For the Shigesada-Kawasaki-Teramoto type population system in \cite{7}, the entropy was
\begin{equation}\label{sktentropy}
\begin{aligned}
h(u)=\sum_{i=1}^{n}\pi_iu_i(\log u_i-1),\quad\pi_i>0.
\end{aligned}
\end{equation}

For the two species chemotaxis system, let us denote $u_3=1-u_1-u_2$, and consider the entropy
\begin{equation}\label{h2u}
\begin{aligned}
h(u)=u_1(\log u_1-1)+u_2(\log u_2-1)+u_3(\log u_3-1).
\end{aligned}
\end{equation}

We remark that $h(u)$ in (\ref{h2u}) is a major functional in showing weak solutions of volume-filling cross-diffusion systems exist. For $u\in\mathcal{D}$, we have
\begin{equation}\label{huchem}
\begin{aligned}
&h^{\prime\prime}(u)A(u)=\Big\{\left(\begin{array}{cccccc}
\frac{1}{u_1} & 0 \\
 0 & \frac{1}{u_2} \\
\end{array}\right)+\frac{1}{u_3}\left(\begin{array}{cccccc}
1 & 1 \\
1 & 1 \\
\end{array}\right)\Big\}\cdot\left(\begin{array}{cccccc}
\alpha_1(1-u_2)+u_2 & (\alpha_1-1)u_1 \\
 (\alpha_2-1)u_2& \alpha_2(1-u_1)+u_1 \\
\end{array}\right)=\\&\left(\begin{array}{cccccc}
\frac{\alpha_1}{u_1}+(1-\alpha_1)\cdot\frac{u_2}{u_1} & \alpha_1-1 \\
 \alpha_2-1 & \frac{\alpha_2}{u_2}+(1-\alpha_2)\cdot\frac{u_1}{u_2} \\
\end{array}\right)+\frac{1}{u_3}\left(\begin{array}{cccccc}
\alpha_1+(\alpha_2-\alpha_1)u_2 & \alpha_2+(\alpha_1-\alpha_2)u_1 \\
\alpha_1+(\alpha_2-\alpha_1)u_2 & \alpha_2+(\alpha_1-\alpha_2)u_1 \\
\end{array}\right).
\end{aligned}
\end{equation}

For general $\alpha_1\neq\alpha_2$ case, we cannot verify assumption (H2) for $h(u)$ and $A(u)$ in (\ref{huchem}). The major novelty of this work is the matrix analysis in Section 2. For each $\alpha_1,\alpha_2$ case, we verify assumptions (H1)-(H3). Then we apply the entropy method to show a weak solution of chemotaxis system exists.

In order to show the weak solution of a cross-diffusion system is unique, we apply the distance functional to measure the ``distance" between two weak solutions of the system. Once the ``distance" between any two weak solutions of the system is zero, then we deduce that the weak solution is unique. 

The distance functional technique was applied in \cite{4,12,13,14}. In \cite{12}, the uniqueness result was shown for a weak solution that is not bounded in general. In \cite{13}, uniqueness results for several typical cross-diffusion systems were discussed.

In general, it is very difficult to construct an appropriate distance functional to show the weak solution of a cross-diffusion system is unique. In \cite{15,16}, weaker uniqueness results for cross-diffusion systems were shown. We call such kind of uniqueness as ``weak-strong uniqueness". 

In these references, so-called relative entropy functionals were applied to measure distances between weak and strong solutions.
In this manuscript, we still focus on the uniqueness proof for weak solutions. 

We assume that $v=(v_1,v_2)$ is also a weak solution of (\ref{equation})-(\ref{initial}). For every test function $\phi=(\phi_1,\phi_2)$, $\phi_1,\phi_2\in L^2(0,T;H^1(\Omega))$, we have
\begin{equation}\label{chemdefv}
\begin{aligned}
&\int_{0}^{T}\langle\partial_tv_1,\phi_1\rangle dt+\int_{Q_T}(\alpha_1\nabla v_1+(\alpha_1-1)(v_1\nabla v_2-v_2\nabla v_1))\cdot\nabla\phi_1dxdt=0,\\&\int_{0}^{T}\langle\partial_tv_2,\phi_2\rangle dt+\int_{Q_T}(\alpha_2\nabla v_2-(\alpha_2-1)(v_1\nabla v_2-v_2\nabla v_1))\cdot\nabla\phi_2dxdt=0,
\end{aligned}
\end{equation}
with the property that $v_i\in L^2(0,T;H^1(\Omega))$, $\partial_tv_i\in L^2(0,T;H^1(\Omega)^{\prime})$, $i=1,2$. The initial datum $v^0=(v^0_1,v^0_2)=u^0$.

Let us take sums of equations in (\ref{chemdef}) and (\ref{chemdefv}). For every test function $\varphi\in L^2(0,T;H^1(\Omega))$, we have
\begin{equation}\label{chementropy}
\begin{aligned}
&\langle\partial_t((\alpha_2-1)u_1+(\alpha_1-1)u_2),\varphi\rangle=\langle\Delta\big(\alpha_1(\alpha_2-1)u_1+\alpha_2(\alpha_1-1)u_2\big),\varphi\rangle,
\\&\langle\partial_t((\alpha_2-1)v_1+(\alpha_1-1)v_2),\varphi\rangle=\langle\Delta\big(\alpha_1(\alpha_2-1)v_1+\alpha_2(\alpha_1-1)v_2\big),\varphi\rangle.
\end{aligned}
\end{equation}
 
In (\ref{chementropy}), we notice that once $\alpha_1=\alpha_2=\alpha$, then we have
\begin{equation}\label{equalchem}
\begin{aligned}
\langle\partial_t(u_1+u_2),\varphi\rangle=\alpha\langle\Delta(u_1+u_2),\varphi\rangle,\quad\langle\partial_t(v_1+v_2),\varphi\rangle=\alpha\langle\Delta(v_1+v_2),\varphi\rangle,\quad u^0=v^0\in\mathcal{D}.
\end{aligned}
\end{equation}

We can show that $u_1+u_2=v_1+v_2$ for a.e. $(x,t)\in\Omega\times(0,T)$ merely by the relation (\ref{equalchem}). For general $\alpha_1\neq\alpha_2$ case, we cannot show that $(\alpha_2-1)u_1+(\alpha_1-1)u_2=(\alpha_2-1)v_1+(\alpha_1-1)v_2$ for a.e. $(x,t)\in\Omega\times(0,T)$ merely by (\ref{chementropy}).

We conjecture that for every $\alpha_1,\alpha_2>0,\alpha_1,\alpha_2\neq 1$, the relation $(\alpha_2-1)u_1+(\alpha_1-1)u_2=(\alpha_2-1)v_1+(\alpha_1-1)v_2$ holds for a.e. $(x,t)\in\Omega\times(0,T)$. It remains one of major challenges in the study of parabolic partial differential equations. Some further remarks will be provided in the last conclusion section.

In \cite{4,12}, the distance functional was chosen as
\begin{equation}\label{oridistance}
\begin{aligned}
d_{\varepsilon}(u,v)=\sum_{i=1,2}\int_{\Omega}\big(\xi_{\varepsilon}(u_i)+\xi_{\varepsilon}(v_i)-2\xi_{\varepsilon}(\frac{u_i+v_i}{2})\big)dx,\quad 0<\varepsilon<1,
\end{aligned}
\end{equation}
where
\begin{equation}\nonumber
\begin{aligned}
\xi_{\varepsilon}(s)=(s+\varepsilon)(\log(s+\varepsilon)-1)+1,\quad s\geq 0.
\end{aligned}
\end{equation}

By the convexity of $\xi_{\varepsilon}$, we have $\xi_{\varepsilon}(u_i)+\xi_{\varepsilon}(v_i)-2\xi_{\varepsilon}(\frac{u_i+v_i}{2})\geq 0$ in $\Omega$, thus $d_{\varepsilon}(u,v)\geq 0$. Since $u_i,v_i$, $i=1,2$ are only nonnegative, we have to add the regularizer $0<\varepsilon<1$ in (\ref{oridistance}). 

In Section 4, we apply the distance functional in (\ref{oridistance}) to show the weak solution of (\ref{equation})-(\ref{initial}) is unique. Another major novelty of this work is reflected in the computation process of the uniqueness proof.

\subsection{Main Steps} 

We divide the existence proof into four steps.

\textbf{Step 1.} In Section 2, we define entropies $h(u)$ in (\ref{entropyepsilon})-(\ref{entropyepsilon4}). We verify that for each $\alpha_1,\alpha_2$ case, $h(u),A(u)$ satisfy assumptions (H1)-(H3). 

The verification of assumption (H2) is reflected in Lemmas \ref{positivesemidefinite}-\ref{positivesemidefinitefourth}, which is the major step of this existence proof.

\textbf{Step 2.} In Section 3, we show a weak solution of (\ref{equation})-(\ref{initial}) exists. We rely on (H1)-(H2) to derive a sequence of approximated solutions indexed by $\tau,\varepsilon>0$. 

In this step, we mainly refer to Lemma 5, \cite{2}. This lemma showed the existence of approximated solutions by applying Lax-Milgram lemma and fixed point theorem.

The proof of Lemma 5, \cite{2} takes several pages and this classical result has been widely applied in global existence analysis of cross-diffusion systems. For these reasons, we briefly state its proof in appendix section.

\textbf{Step 3.} We provide uniform estimates of approximated solutions indexed by $\tau,\varepsilon>0$. By (H2)-(H3), we derive estimation formulas (\ref{firstaubin})-(\ref{secondaubin}).

\textbf{Step 4.} We let $(\tau,\varepsilon)\to 0$. We apply Aubin-Lions lemma in \cite{17,18} to show that the sequence of approximated solutions converges to a weak solution of (\ref{equation})-(\ref{initial}).

Then we divide the uniqueness proof into these three steps.

\textbf{Step 1.} In Section 4, we list Lemma \ref{j2require} as an important preparation. We show that once $\alpha_1=\alpha_2=\alpha$, then $u_1+u_2=v_1+v_2$, a.e. $(x,t)\in\Omega\times(0,T)$. We for convenience denote $u_3=1-u_1-u_2$, $v_3=1-v_1-v_2$.

\textbf{Step 2.} We apply the distance functional in (\ref{oridistance}) to measure the ``distance" between $u=(u_1,u_2)$ and $v=(v_1,v_2)$. $u$ and $v$ are two weak solutions of (\ref{equation})-(\ref{initial}).

\textbf{Step 3.} The distance functional in (\ref{oridistance}) is indexed by $0<\varepsilon<1$. We let $\varepsilon\to 0$, to show that $u=v$, a.e. $(x,t)\in\Omega\times(0,T)$.

\section{Matrix Analysis}

\subsection{The Construction of $h(u)$} Let us consider the case when $0<\alpha_1,\alpha_2<1$. We define the entropy $h(u)$ as
\begin{equation}\label{entropyepsilon}
\begin{aligned}
h(u)=\frac{\beta u_1(\log u_1-1)}{1-\alpha_1}+\frac{\beta u_2(\log u_2-1)}{1-\alpha_2}+u_3(\log u_3-1)+\frac{\beta}{1-\alpha_1}+\frac{\beta}{1-\alpha_2}+1,\quad\beta>0.
\end{aligned}
\end{equation}

We verify that $h\in C^2(\mathcal{D},(0,\infty))$, and
\begin{equation}\nonumber
\begin{aligned}
\partial h/\partial u_i=\beta\log u_i/(1-\alpha_i)-\log u_3,\quad i=1,2.
\end{aligned}
\end{equation}

For every $(w_1,w_2)\in\mathbb{R}^2$, we define 
\begin{equation}\nonumber
\begin{aligned}
g(y)=e^{(1-\alpha_1)w_1/\beta}(1-y)^{(1-\alpha_1)/\beta}+e^{(1-\alpha_2)w_2/\beta}(1-y)^{(1-\alpha_2)/\beta},\quad 0<y<1.
\end{aligned}
\end{equation}
$g$ is continuous and nonincreasing in $0<y<1$, and
\begin{equation}\nonumber
\begin{aligned}
g(0)=e^{(1-\alpha_1)w_1/\beta}+e^{(1-\alpha_2)w_2/\beta}>0,\quad g(1)=0.
\end{aligned}
\end{equation}
We deduce that there exists a unique fixed point $0<y_0<1$, such that $g(y_0)=y_0$. 

Then we define $u_i=e^{(1-\alpha_i)w_i/\beta}(1-y_0)^{(1-\alpha_i)/\beta}>0$($i=1,2$), which satisfies $u_1+u_2=g(y_0)=y_0<1$. Consequently, $u=(u_1,u_2)\in\mathcal{D}$, $u_3=1-u_1-u_2=1-y_0$, and $$w_i=\beta\log u_i/(1-\alpha_i)-\log u_3=\partial h/\partial u_i,\quad i=1,2.$$ Then we show that $h^{\prime}:\mathcal{D}\to\mathbb{R}^2$ is invertible on $\mathbb{R}^2$, which satisfies assumption (H1). 

Entropies in (\ref{entropyepsilon2})-(\ref{entropyepsilon4}) also satisfy assumption (H1), with similar verification processes. In Lemmas \ref{positivesemidefinite}-\ref{positivesemidefinitefourth}, we verify assumption (H2). We divide the verification into four cases.

\subsection{Case 1: $0<\alpha_1,\alpha_2<1$}

Let us verify assumption (H2) in Lemma \ref{positivesemidefinite}.

\begin{lem}\label{positivesemidefinite}
For every $0<\alpha_1,\alpha_2<1$, $z=(z_1,z_2)\in\mathbb{R}^2$, $u\in\mathcal{D}$, $h(u)$ given by $(\ref{entropyepsilon})$, we have
\begin{equation}\nonumber
\begin{aligned}
z^{\mathrm{T}}h^{\prime\prime}(u)A(u)z\geq c(z^2_1+z^2_2),
\end{aligned}
\end{equation}
with $c>0$ a constant independent of $u$ and $z$.
\end{lem}

\begin{proof}

By (\ref{entropyepsilon}), we have
\begin{equation}\nonumber
\begin{aligned}
&h^{\prime\prime}(u)=\left(\begin{array}{cccccc}
\frac{\beta}{(1-\alpha_1)u_1} & 0 \\
 0 & \frac{\beta}{(1-\alpha_2)u_2} \\
\end{array}\right)+\frac{1}{u_3}\left(\begin{array}{cccccc}
1 & 1 \\
1 & 1 \\
\end{array}\right).
\end{aligned}
\end{equation}

In the following computation, once $\alpha_1<\alpha_2$, we decompose the matrix as
\begin{equation}\nonumber
\begin{aligned}
&\frac{1}{u_3}\left(\begin{array}{cccccc}
\alpha_1+(\alpha_2-\alpha_1)u_2 & \alpha_2+(\alpha_1-\alpha_2)u_1 \\
\alpha_1+(\alpha_2-\alpha_1)u_2 & \alpha_2+(\alpha_1-\alpha_2)u_1 \\
\end{array}\right)\\&=\frac{\alpha_1+(\alpha_2-\alpha_1)u_2}{u_3}\left(\begin{array}{cccccc}
1 & 1 \\
1 & 1 \\
\end{array}\right)+(\alpha_2-\alpha_1)\left(\begin{array}{cccccc}
0 & 1 \\
0 & 1 \\
\end{array}\right).
\end{aligned}
\end{equation}

If $\alpha_1>\alpha_2$, the decomposition is given as
\begin{equation}\nonumber
\begin{aligned}
&\frac{1}{u_3}\left(\begin{array}{cccccc}
\alpha_1+(\alpha_2-\alpha_1)u_2 & \alpha_2+(\alpha_1-\alpha_2)u_1 \\
\alpha_1+(\alpha_2-\alpha_1)u_2 & \alpha_2+(\alpha_1-\alpha_2)u_1 \\
\end{array}\right)\\&=\frac{\alpha_2+(\alpha_1-\alpha_2)u_1}{u_3}\left(\begin{array}{cccccc}
1 & 1 \\
1 & 1 \\
\end{array}\right)+(\alpha_1-\alpha_2)\left(\begin{array}{cccccc}
1 & 0 \\
1 & 0 \\
\end{array}\right).
\end{aligned}
\end{equation}

Once $\alpha_1<\alpha_2$, we have
\begin{equation}\nonumber
\begin{aligned}
&h^{\prime\prime}(u)A(u)=\\&\Big\{\left(\begin{array}{cccccc}
\frac{\beta}{(1-\alpha_1)u_1} & 0 \\
 0 & \frac{\beta}{(1-\alpha_2)u_2} \\
\end{array}\right)+\frac{1}{u_3}\left(\begin{array}{cccccc}
1 & 1 \\
1 & 1 \\
\end{array}\right)\Big\}\cdot\left(\begin{array}{cccccc}
\alpha_1+(1-\alpha_1)u_2 & (\alpha_1-1)u_1 \\
(\alpha_2-1)u_2 & \alpha_2+(1-\alpha_2)u_1 \\
\end{array}\right)=\\&\beta\left(\begin{array}{cccccc}
\frac{\alpha_1}{(1-\alpha_1)u_1}+\frac{u_2}{u_1} & -1 \\
-1 & \frac{\alpha_2}{(1-\alpha_2)u_2}+\frac{u_1}{u_2} \\
\end{array}\right)+\frac{\alpha_1+(\alpha_2-\alpha_1)u_2}{u_3}\left(\begin{array}{cccccc}
1 & 1 \\
1 & 1 \\
\end{array}\right)+(\alpha_2-\alpha_1)\left(\begin{array}{cccccc}
0 & 1 \\
0 & 1 \\
\end{array}\right).
\end{aligned}
\end{equation}

Let us choose the constant $\beta$ as
\begin{equation}\nonumber
\begin{aligned}
\beta=\frac{(1-\alpha_1)(\alpha_2-\alpha_1)}{2\alpha_1},\text{ and denote }c=\min\Big\{\frac{\alpha_2-\alpha_1}{4},\frac{(1-\alpha_1)(\alpha_2-\alpha_1)\alpha_2}{2(1-\alpha_2)\alpha_1}\Big\},
\end{aligned}
\end{equation}
which follows that
\begin{equation}\nonumber
\begin{aligned}
&\frac{\beta\alpha_1z^2_1}{(1-\alpha_1)u_1}+\Big(\frac{\beta\alpha_2}{(1-\alpha_2)u_2}+(\alpha_2-\alpha_1)\Big)z^2_2+(\alpha_2-\alpha_1)z_1z_2\\&=(\alpha_2-\alpha_1)(\frac{z_1}{2}+z_2)^2+\Big(\frac{\beta\alpha_1}{(1-\alpha_1)u_1}-\frac{\alpha_2-\alpha_1}{4}\Big)z^2_1+\frac{\beta\alpha_2z^2_2}{(1-\alpha_2)u_2}\\&\geq\frac{(\alpha_2-\alpha_1)z^2_1}{4}+\frac{(1-\alpha_1)(\alpha_2-\alpha_1)\alpha_2z^2_2}{2(1-\alpha_2)\alpha_1}\geq c(z^2_1+z^2_2).
\end{aligned}
\end{equation}

For every $u\in\mathcal{D}$, $z=(z_1,z_2)\in\mathbb{R}^2$, we have
\begin{equation}\nonumber
\begin{aligned}
&z^{\mathrm{T}}h^{\prime\prime}(u)A(u)z=\frac{\beta\alpha_1z^2_1}{(1-\alpha_1)u_1}+\Big(\frac{\beta\alpha_2}{(1-\alpha_2)u_2}+(\alpha_2-\alpha_1)\Big)z^2_2+(\alpha_2-\alpha_1)z_1z_2\\&+\beta\big(\sqrt{\frac{u_2}{u_1}}z_1-\sqrt{\frac{u_1}{u_2}}z_2\big)^2+\frac{\alpha_1+(\alpha_2-\alpha_1)u_2}{u_3}\cdot(z_1+z_2)^2\geq c(z^2_1+z^2_2).
\end{aligned}
\end{equation}

Once $\alpha_1>\alpha_2$, we have
\begin{equation}\nonumber
\begin{aligned}
&h^{\prime\prime}(u)A(u)=\\&\beta\left(\begin{array}{cccccc}
\frac{\alpha_1}{(1-\alpha_1)u_1}+\frac{u_2}{u_1} & -1 \\
-1 & \frac{\alpha_2}{(1-\alpha_2)u_2}+\frac{u_1}{u_2} \\
\end{array}\right)+\frac{\alpha_2+(\alpha_1-\alpha_2)u_1}{u_3}\left(\begin{array}{cccccc}
1 & 1 \\
1 & 1 \\
\end{array}\right)+(\alpha_1-\alpha_2)\left(\begin{array}{cccccc}
1 & 0 \\
1 & 0 \\
\end{array}\right).
\end{aligned}
\end{equation}

Let us choose the constant $\beta$ as
\begin{equation}\nonumber
\begin{aligned}
\beta=\frac{(1-\alpha_2)(\alpha_1-\alpha_2)}{2\alpha_2},\text{ and denote }c=\min\Big\{\frac{\alpha_1-\alpha_2}{4},\frac{(1-\alpha_2)(\alpha_1-\alpha_2)\alpha_1}{2(1-\alpha_1)\alpha_2}\Big\},
\end{aligned}
\end{equation}
which follows that
\begin{equation}\nonumber
\begin{aligned}
&\Big(\frac{\beta\alpha_1}{(1-\alpha_1)u_1}+(\alpha_1-\alpha_2)\Big)z^2_1+\frac{\beta\alpha_2z^2_2}{(1-\alpha_2)u_2}+(\alpha_1-\alpha_2)z_1z_2\\&=(\alpha_1-\alpha_2)(z_1+\frac{z_2}{2})^2+\frac{\beta\alpha_1z^2_1}{(1-\alpha_1)u_1}+\Big(\frac{\beta\alpha_2}{(1-\alpha_2)u_2}-\frac{\alpha_1-\alpha_2}{4}\Big)z^2_2\\&\geq\frac{(1-\alpha_2)(\alpha_1-\alpha_2)\alpha_1z^2_1}{2(1-\alpha_1)\alpha_2}+\frac{(\alpha_1-\alpha_2)z^2_2}{4}\geq c(z^2_1+z^2_2).
\end{aligned}
\end{equation}

For every $u\in\mathcal{D}$, $z=(z_1,z_2)\in\mathbb{R}^2$, we have
\begin{equation}\nonumber
\begin{aligned}
&z^{\mathrm{T}}h^{\prime\prime}(u)A(u)z=\Big(\frac{\beta\alpha_1}{(1-\alpha_1)u_1}+(\alpha_1-\alpha_2)\Big)z^2_1+\frac{\beta\alpha_2z^2_2}{(1-\alpha_2)u_2}+(\alpha_1-\alpha_2)z_1z_2\\&+\beta\big(\sqrt{\frac{u_2}{u_1}}z_1-\sqrt{\frac{u_1}{u_2}}z_2\big)^2+\frac{\alpha_2+(\alpha_1-\alpha_2)u_1}{u_3}\cdot(z_1+z_2)^2\geq c(z^2_1+z^2_2).
\end{aligned}
\end{equation}

Once $\alpha_1=\alpha_2=\alpha$, we choose $\beta=1$, and denote $c=\alpha/(1-\alpha)$. For every $u\in\mathcal{D}$, $z=(z_1,z_2)\in\mathbb{R}^2$, we have
\begin{equation}\nonumber
\begin{aligned}
&z^{\mathrm{T}}h^{\prime\prime}(u)A(u)z=\\&\frac{\alpha z^2_1}{(1-\alpha)u_1}+\frac{\alpha z^2_2}{(1-\alpha)u_2}+\big(\sqrt{\frac{u_2}{u_1}}z_1-\sqrt{\frac{u_1}{u_2}}z_2\big)^2+\frac{\alpha}{u_3}\cdot(z_1+z_2)^2\geq c(z^2_1+z^2_2).
\end{aligned}
\end{equation}

Then we finish the proof of Lemma \ref{positivesemidefinite}.
\end{proof}

\subsection{Case 2: $\alpha_1,\alpha_2>1$} Let us define the entropy $h(u)$ as
\begin{equation}\label{entropyepsilon2}
\begin{aligned}
h(u)=\frac{\beta u_1(\log u_1-1)}{\alpha_1-1}+\frac{\beta u_2(\log u_2-1)}{\alpha_2-1}+u_3(\log u_3-1)+\frac{\beta}{\alpha_1-1}+\frac{\beta}{\alpha_2-1}+1,\quad\beta>0,
\end{aligned}
\end{equation}
when $\alpha_1,\alpha_2>1$. Then we show Lemma \ref{positivesemidefinitesecond}.

\begin{lem}\label{positivesemidefinitesecond}
For every $\alpha_1,\alpha_2>1$, $z=(z_1,z_2)\in\mathbb{R}^2$, $u\in\mathcal{D}$, $h(u)$ given by $(\ref{entropyepsilon2})$, we have
\begin{equation}\nonumber
\begin{aligned}
z^{\mathrm{T}}h^{\prime\prime}(u)A(u)z\geq c(z^2_1+z^2_2),
\end{aligned}
\end{equation}
with $c>0$ a constant independent of $u$ and $z$.
\end{lem}

\begin{proof}

By (\ref{entropyepsilon2}), we have
\begin{equation}\nonumber
\begin{aligned}
&h^{\prime\prime}(u)=\left(\begin{array}{cccccc}
\frac{\beta}{(\alpha_1-1)u_1} & 0 \\
 0 & \frac{\beta}{(\alpha_2-1)u_2} \\
\end{array}\right)+\frac{1}{u_3}\left(\begin{array}{cccccc}
1 & 1 \\
1 & 1 \\
\end{array}\right).
\end{aligned}
\end{equation}

Once $\alpha_1<\alpha_2$, we have
\begin{equation}\nonumber
\begin{aligned}
&h^{\prime\prime}(u)A(u)=\\&\Big\{\left(\begin{array}{cccccc}
\frac{\beta}{(\alpha_1-1)u_1} & 0 \\
0 & \frac{\beta}{(\alpha_2-1)u_2} \\
\end{array}\right)+\frac{1}{u_3}\left(\begin{array}{cccccc}
1 & 1 \\
1 & 1 \\
\end{array}\right)\Big\}\cdot\left(\begin{array}{cccccc}
\alpha_1+(1-\alpha_1)u_2 & (\alpha_1-1)u_1 \\
(\alpha_2-1)u_2 & \alpha_2+(1-\alpha_2)u_1 \\
\end{array}\right)\\&=\beta\left(\begin{array}{cccccc}
\frac{1}{(\alpha_1-1)u_1}+\frac{u_3}{u_1}+1 & 1 \\
1 & \frac{1}{(\alpha_2-1)u_2}+\frac{u_3}{u_2}+1 \\
\end{array}\right)+\frac{\alpha_1+(\alpha_2-\alpha_1)u_2}{u_3}\left(\begin{array}{cccccc}
1 & 1 \\
1 & 1 \\
\end{array}\right)\\&+(\alpha_2-\alpha_1)\left(\begin{array}{cccccc}
0 & 1 \\
0 & 1 \\
\end{array}\right).
\end{aligned}
\end{equation}

Let us choose the constant $\beta$ as
\begin{equation}\nonumber
\begin{aligned}
\beta=\frac{(\alpha_1-1)(\alpha_2-\alpha_1)}{2},\text{ and denote }c=\min\Big\{\frac{\alpha_2-\alpha_1}{4},\frac{(\alpha_1-1)(\alpha_2-\alpha_1)}{2(\alpha_2-1)}\Big\},
\end{aligned}
\end{equation}
which follows that
\begin{equation}\nonumber
\begin{aligned}
&\frac{\beta z^2_1}{(\alpha_1-1)u_1}+\Big(\frac{\beta}{(\alpha_2-1)u_2}+(\alpha_2-\alpha_1)\Big)z^2_2+(\alpha_2-\alpha_1)z_1z_2\\&=(\alpha_2-\alpha_1)(\frac{z_1}{2}+z_2)^2+\Big(\frac{\beta}{(\alpha_1-1)u_1}-\frac{\alpha_2-\alpha_1}{4}\Big)z^2_1+\frac{\beta z^2_2}{(\alpha_2-1)u_2}\\&\geq\frac{(\alpha_2-\alpha_1)z^2_1}{4}+\frac{(\alpha_1-1)(\alpha_2-\alpha_1)z^2_2}{2(\alpha_2-1)}\geq c(z^2_1+z^2_2).
\end{aligned}
\end{equation}

For every $u\in\mathcal{D}$, $z=(z_1,z_2)\in\mathbb{R}^2$, we have
\begin{equation}\nonumber
\begin{aligned}
&z^{\mathrm{T}}h^{\prime\prime}(u)A(u)z=\frac{\beta z^2_1}{(\alpha_1-1)u_1}+\Big(\frac{\beta}{(\alpha_2-1)u_2}+(\alpha_2-\alpha_1)\Big)z^2_2+(\alpha_2-\alpha_1)z_1z_2\\&+\beta u_3\big(\frac{z^2_1}{u_1}+\frac{z^2_2}{u_2}\big)+\big(\beta+\frac{\alpha_1+(\alpha_2-\alpha_1)u_2}{u_3}\big)\cdot(z_1+z_2)^2\geq c(z^2_1+z^2_2).
\end{aligned}
\end{equation}

Once $\alpha_1>\alpha_2$, we have
\begin{equation}\nonumber
\begin{aligned}
&h^{\prime\prime}(u)A(u)=\beta\left(\begin{array}{cccccc}
\frac{1}{(\alpha_1-1)u_1}+\frac{u_3}{u_1}+1 & 1 \\
1 & \frac{1}{(\alpha_2-1)u_2}+\frac{u_3}{u_2}+1 \\
\end{array}\right)\\&+\frac{\alpha_2+(\alpha_1-\alpha_2)u_1}{u_3}\left(\begin{array}{cccccc}
1 & 1 \\
1 & 1 \\
\end{array}\right)+(\alpha_1-\alpha_2)\left(\begin{array}{cccccc}
1 & 0 \\
1 & 0 \\
\end{array}\right).
\end{aligned}
\end{equation}

Let us choose the constant $\beta$ as
\begin{equation}\nonumber
\begin{aligned}
\beta=\frac{(\alpha_2-1)(\alpha_1-\alpha_2)}{2},\text{ and denote }c=\min\Big\{\frac{\alpha_1-\alpha_2}{4},\frac{(\alpha_2-1)(\alpha_1-\alpha_2)}{2(\alpha_1-1)}\Big\},
\end{aligned}
\end{equation}
which follows that
\begin{equation}\nonumber
\begin{aligned}
&\Big(\frac{\beta}{(\alpha_1-1)u_1}+(\alpha_1-\alpha_2)\Big)z^2_1+\frac{\beta z^2_2}{(\alpha_2-1)u_2}+(\alpha_1-\alpha_2)z_1z_2\\&=(\alpha_1-\alpha_2)(z_1+\frac{z_2}{2})^2+\frac{\beta z^2_1}{(\alpha_1-1)u_1}+\Big(\frac{\beta}{(\alpha_2-1)u_2}-\frac{\alpha_1-\alpha_2}{4}\Big)z^2_2\\&\geq\frac{(\alpha_2-1)(\alpha_1-\alpha_2)z^2_1}{2(\alpha_1-1)}+\frac{(\alpha_1-\alpha_2)z^2_2}{4}\geq c(z^2_1+z^2_2).
\end{aligned}
\end{equation}

For every $u\in\mathcal{D}$, $z=(z_1,z_2)\in\mathbb{R}^2$, we have
\begin{equation}\nonumber
\begin{aligned}
&z^{\mathrm{T}}h^{\prime\prime}(u)A(u)z=\Big(\frac{\beta}{(\alpha_1-1)u_1}+(\alpha_1-\alpha_2)\Big)z^2_1+\frac{\beta z^2_2}{(\alpha_2-1)u_2}+(\alpha_1-\alpha_2)z_1z_2\\&+\beta u_3\big(\frac{z^2_1}{u_1}+\frac{z^2_2}{u_2}\big)+\big(\beta+\frac{\alpha_2+(\alpha_1-\alpha_2)u_1}{u_3}\big)\cdot(z_1+z_2)^2\geq c(z^2_1+z^2_2).
\end{aligned}
\end{equation}

Once $\alpha_1=\alpha_2=\alpha$, we choose $\beta=1$, and denote $c=1/(\alpha-1)$. For every $u\in\mathcal{D}$, $z=(z_1,z_2)\in\mathbb{R}^2$, we have
\begin{equation}\nonumber
\begin{aligned}
&z^{\mathrm{T}}h^{\prime\prime}(u)A(u)z=\\&\frac{z^2_1}{(\alpha-1)u_1}+\frac{z^2_2}{(\alpha-1)u_2}+\big(\frac{\alpha}{u_3}+1\big)\cdot(z_1+z_2)^2+u_3\big(\frac{z^2_1}{u_1}+\frac{z^2_2}{u_2}\big)\geq c(z^2_1+z^2_2).
\end{aligned}
\end{equation}

Then we have shown Lemma \ref{positivesemidefinitesecond}.

\end{proof}

\subsection{Case 3: $0<\alpha_2<1<\alpha_1$.} If $0<\alpha_2<1<\alpha_1$, we define the entropy $h(u)$ as
\begin{equation}\label{entropyepsilon3}
\begin{aligned}
h(u)=\frac{\beta u_1(\log u_1-1)}{\alpha_1-1}+\frac{\beta u_2(\log u_2-1)}{1-\alpha_2}+u_3(\log u_3-1)+\frac{\beta}{\alpha_1-1}+\frac{\beta}{1-\alpha_2}+1,\quad\beta>0.
\end{aligned}
\end{equation}

\begin{lem}\label{positivesemidefinitethird}
For every $0<\alpha_2<1<\alpha_1$, $z=(z_1,z_2)\in\mathbb{R}^2$, $u\in\mathcal{D}$, $h(u)$ given by $(\ref{entropyepsilon3})$, we have
\begin{equation}\nonumber
\begin{aligned}
z^{\mathrm{T}}h^{\prime\prime}(u)A(u)z\geq c(z^2_1+z^2_2),
\end{aligned}
\end{equation}
with $c>0$ a constant independent of $u$ and $z$.
\end{lem}

\begin{proof}

By (\ref{entropyepsilon3}), we have
\begin{equation}\nonumber
\begin{aligned}
&h^{\prime\prime}(u)=\left(\begin{array}{cccccc}
\frac{\beta}{(\alpha_1-1)u_1} & 0 \\
 0 & \frac{\beta}{(1-\alpha_2)u_2} \\
\end{array}\right)+\frac{1}{u_3}\left(\begin{array}{cccccc}
1 & 1 \\
1 & 1 \\
\end{array}\right),
\end{aligned}
\end{equation}
which follows
\begin{equation}\nonumber
\begin{aligned}
&h^{\prime\prime}(u)A(u)=\\&\Big\{\left(\begin{array}{cccccc}
\frac{\beta}{(\alpha_1-1)u_1} & 0 \\
0 & \frac{\beta}{(1-\alpha_2)u_2} \\
\end{array}\right)+\frac{1}{u_3}\left(\begin{array}{cccccc}
1 & 1 \\
1 & 1 \\
\end{array}\right)\Big\}\cdot\left(\begin{array}{cccccc}
\alpha_1+(1-\alpha_1)u_2 & (\alpha_1-1)u_1 \\
(\alpha_2-1)u_2 & \alpha_2+(1-\alpha_2)u_1 \\
\end{array}\right)\\&=\beta\left(\begin{array}{cccccc}
\frac{1}{(\alpha_1-1)u_1}+\frac{u_3}{u_1}+1 & 1 \\
-1 & \frac{\alpha_2}{(1-\alpha_2)u_2}+\frac{u_1}{u_2} \\
\end{array}\right)+\frac{\alpha_2+(\alpha_1-\alpha_2)u_1}{u_3}\left(\begin{array}{cccccc}
1 & 1 \\
1 & 1 \\
\end{array}\right)\\&+(\alpha_1-\alpha_2)\left(\begin{array}{cccccc}
1 & 0 \\
1 & 0 \\
\end{array}\right).
\end{aligned}
\end{equation}

Let us choose the constant $\beta$ as
\begin{equation}\nonumber
\begin{aligned}
\beta=\frac{(1-\alpha_2)(\alpha_1-\alpha_2)}{2\alpha_2},\text{ and denote }c=\min\Big\{\frac{\alpha_1-\alpha_2}{4},\frac{(1-\alpha_2)(\alpha_1-\alpha_2)\alpha_1}{2(\alpha_1-1)\alpha_2}\Big\},
\end{aligned}
\end{equation}
which follows that
\begin{equation}\nonumber
\begin{aligned}
&\Big(\frac{\beta}{(\alpha_1-1)u_1}+(\alpha_1-\alpha_2+\beta)\Big)z^2_1+\frac{\beta\alpha_2z^2_2}{(1-\alpha_2)u_2}+(\alpha_1-\alpha_2)z_1z_2\\&=(\alpha_1-\alpha_2)(z_1+\frac{z_2}{2})^2+\beta z^2_1+\frac{\beta z^2_1}{(\alpha_1-1)u_1}+\Big(\frac{\beta\alpha_2}{(1-\alpha_2)u_2}-\frac{\alpha_1-\alpha_2}{4}\Big)z^2_2\\&\geq\frac{(1-\alpha_2)(\alpha_1-\alpha_2)\alpha_1z^2_1}{2(\alpha_1-1)\alpha_2}+\frac{(\alpha_1-\alpha_2)z^2_2}{4}\geq c(z^2_1+z^2_2).
\end{aligned}
\end{equation}

For every $u\in\mathcal{D}$, $z=(z_1,z_2)\in\mathbb{R}^2$, we have
\begin{equation}\nonumber
\begin{aligned}
&z^{\mathrm{T}}h^{\prime\prime}(u)A(u)z=\Big(\frac{\beta}{(\alpha_1-1)u_1}+(\alpha_1-\alpha_2+\beta)\Big)z^2_1+\frac{\beta\alpha_2z^2_2}{(1-\alpha_2)u_2}+(\alpha_1-\alpha_2)z_1z_2\\&+\frac{\beta u_3}{u_1}\cdot z^2_1+\frac{\beta u_1}{u_2}\cdot z^2_2+\frac{\alpha_2+(\alpha_1-\alpha_2)u_1}{u_3}\cdot(z_1+z_2)^2\geq c(z^2_1+z^2_2),
\end{aligned}
\end{equation}
then we have shown Lemma \ref{positivesemidefinitethird}.

\end{proof}

\subsection{Case 4: $0<\alpha_1<1<\alpha_2$.} If $0<\alpha_1<1<\alpha_2$, we define the entropy $h(u)$ as
\begin{equation}\label{entropyepsilon4}
\begin{aligned}
h(u)=\frac{\beta u_1(\log u_1-1)}{1-\alpha_1}+\frac{\beta u_2(\log u_2-1)}{\alpha_2-1}+u_3(\log u_3-1)+\frac{\beta}{1-\alpha_1}+\frac{\beta}{\alpha_2-1}+1,\quad\beta>0.
\end{aligned}
\end{equation}

\begin{lem}\label{positivesemidefinitefourth}
For every $0<\alpha_1<1<\alpha_2$, $z=(z_1,z_2)\in\mathbb{R}^2$, $u\in\mathcal{D}$, $h(u)$ given by $(\ref{entropyepsilon4})$, we have
\begin{equation}\nonumber
\begin{aligned}
z^{\mathrm{T}}h^{\prime\prime}(u)A(u)z\geq c(z^2_1+z^2_2),
\end{aligned}
\end{equation}
with $c>0$ a constant independent of $u$ and $z$.
\end{lem}

\begin{proof}

By (\ref{entropyepsilon4}), we have
\begin{equation}\nonumber
\begin{aligned}
&h^{\prime\prime}(u)=\left(\begin{array}{cccccc}
\frac{\beta}{(1-\alpha_1)u_1} & 0 \\
 0 & \frac{\beta}{(\alpha_2-1)u_2} \\
\end{array}\right)+\frac{1}{u_3}\left(\begin{array}{cccccc}
1 & 1 \\
1 & 1 \\
\end{array}\right),
\end{aligned}
\end{equation}
which follows
\begin{equation}\nonumber
\begin{aligned}
&h^{\prime\prime}(u)A(u)=\\&\Big\{\left(\begin{array}{cccccc}
\frac{\beta}{(1-\alpha_1)u_1} & 0 \\
0 & \frac{\beta}{(\alpha_2-1)u_2} \\
\end{array}\right)+\frac{1}{u_3}\left(\begin{array}{cccccc}
1 & 1 \\
1 & 1 \\
\end{array}\right)\Big\}\cdot\left(\begin{array}{cccccc}
\alpha_1+(1-\alpha_1)u_2 & (\alpha_1-1)u_1 \\
(\alpha_2-1)u_2 & \alpha_2+(1-\alpha_2)u_1 \\
\end{array}\right)\\&=\beta\left(\begin{array}{cccccc}
\frac{\alpha_1}{(1-\alpha_1)u_1}+\frac{u_2}{u_1} & -1 \\
1 & \frac{1}{(\alpha_2-1)u_2}+\frac{u_3}{u_2}+1 \\
\end{array}\right)+\frac{\alpha_1+(\alpha_2-\alpha_1)u_2}{u_3}\left(\begin{array}{cccccc}
1 & 1 \\
1 & 1 \\
\end{array}\right)\\&+(\alpha_2-\alpha_1)\left(\begin{array}{cccccc}
0 & 1 \\
0 & 1 \\
\end{array}\right).
\end{aligned}
\end{equation}

Let us choose the constant $\beta$ as
\begin{equation}\nonumber
\begin{aligned}
\beta=\frac{(1-\alpha_1)(\alpha_2-\alpha_1)}{2\alpha_1},\text{ and denote }c=\min\Big\{\frac{\alpha_2-\alpha_1}{4},\frac{(1-\alpha_1)(\alpha_2-\alpha_1)\alpha_2}{2(\alpha_2-1)\alpha_1}\Big\},
\end{aligned}
\end{equation}
which follows that
\begin{equation}\nonumber
\begin{aligned}
&\frac{\beta\alpha_1z^2_1}{(1-\alpha_1)u_1}+\Big(\frac{\beta}{(\alpha_2-1)u_2}+(\alpha_2-\alpha_1+\beta)\Big)z^2_2+(\alpha_2-\alpha_1)z_1z_2\\&=(\alpha_2-\alpha_1)(\frac{z_1}{2}+z_2)^2+\Big(\frac{\beta\alpha_1}{(1-\alpha_1)u_1}-\frac{\alpha_2-\alpha_1}{4}\Big)z^2_1+\beta z^2_2+\frac{\beta z^2_2}{(\alpha_2-1)u_2}\\&\geq\frac{(\alpha_2-\alpha_1)z^2_1}{4}+\frac{(1-\alpha_1)(\alpha_2-\alpha_1)\alpha_2z^2_2}{2(\alpha_2-1)\alpha_1}\geq c(z^2_1+z^2_2).
\end{aligned}
\end{equation}

For every $u\in\mathcal{D}$, $z=(z_1,z_2)\in\mathbb{R}^2$, we have
\begin{equation}\nonumber
\begin{aligned}
&z^{\mathrm{T}}h^{\prime\prime}(u)A(u)z=\frac{\beta\alpha_1z^2_1}{(1-\alpha_1)u_1}+\Big(\frac{\beta}{(\alpha_2-1)u_2}+(\alpha_2-\alpha_1+\beta)\Big)z^2_2+(\alpha_2-\alpha_1)z_1z_2\\&+\frac{\beta u_2}{u_1}\cdot z^2_1+\frac{\beta u_3}{u_2}\cdot z^2_2+\frac{\alpha_1+(\alpha_2-\alpha_1)u_2}{u_3}\cdot(z_1+z_2)^2\geq c(z^2_1+z^2_2),
\end{aligned}
\end{equation}
then we have shown Lemma \ref{positivesemidefinitefourth}.

\end{proof}

\section{Proof of Theorem \ref{mainresult}}

Let us give the proof of Theorem \ref{mainresult}.

\begin{proof}

\textbf{Step 1.} Let $T>0$, $N\in\mathbb{N}$, $\tau=T/N$, $0<\varepsilon<1$, and $m\in\mathbb{N}$ with $m>d/2$. This ensures that the embedding $H^m(\Omega)\hookrightarrow L^{\infty}(\Omega)$ is compact. 

Lemmas \ref{positivesemidefinite}-\ref{positivesemidefinitefourth} indicate that assumption (H2) can be satisfied. By Lemma 5 of \cite{2}, we can show that given $w^{k-1}\in L^{\infty}(\Omega)$ for $k\in\mathbb{N}$, there exists $w^k\in H^m(\Omega)$, $u(w^k(x))\in\mathcal{D}$ for $x\in\Omega$, such that
\begin{equation}\label{inductionhyp}
\begin{aligned}
&\frac{1}{\tau}\int_{\Omega}(u(w^k)-u(w^{k-1}))\cdot\phi dx+\int_{\Omega}\nabla\phi:B(w^k)\nabla w^kdx\\&+\varepsilon\int_{\Omega}\Big(\sum_{|\alpha|=m}D^{\alpha}w^k\cdot D^{\alpha}\phi+w^k\cdot\phi\Big)dx=0,
\end{aligned}
\end{equation}
for every $\phi\in H^m(\Omega)$. Here, $A:B=\sum_{i,j}a_{ij}b_{ij}$ for matrices $A=(a_{ij})$ and $B=(b_{ij})$.

In (\ref{inductionhyp}), $u(w^k)=(h^{\prime})^{-1}(w^k)$, $B(w^k)=A(u(w^k))H(u(w^k))^{-1}$, and $H(u)=h^{\prime\prime}(u)$. $|\alpha|=\alpha_1+\cdot\cdot\cdot+\alpha_d=m$ is a multiindex and $D^{\alpha}=\partial^{|\alpha|}/(\partial x^{\alpha_1}_1\cdot\cdot\cdot\partial x^{\alpha_d}_d)$ is a partial derivative of order $m$. If $k=1$, we define $w^0=h^{\prime}(u^0)$.

Let us denote $u^k=u(w^k)$, and introduce the piecewise in time constant functions $u^{(\tau)}=u^k$, $w^{(\tau)}=w^k$, for $x\in\Omega$, $t\in((k-1)\tau,k\tau]$. At time $t=0$, we set $u^{(\tau)}(\cdot,0)=(h^{\prime})^{-1}(u^0)$. Let $u^{(\tau)}=(u^{(\tau)}_1,u^{(\tau)}_2)$, we define the backward shift operator $(\sigma_{\tau}u^{(\tau)})(x,t)=u^{k-1}$ for $x\in\Omega$, $t\in((k-1)\tau,k\tau]$. $u^{(\tau)}$ satisfies that
\begin{equation}\label{inductiontau}
\begin{aligned}
&\frac{1}{\tau}\int_{0}^{T}\int_{\Omega}(u^{(\tau)}-\sigma_{\tau}u^{(\tau)})\cdot\phi dxdt+\int_{0}^{T}\int_{\Omega}\nabla\phi:B(w^{(\tau)})\nabla w^{(\tau)}dxdt\\&+\varepsilon\int_{0}^{T}\int_{\Omega}\Big(\sum_{|\alpha|=m}D^{\alpha}w^{(\tau)}\cdot D^{\alpha}\phi+w^{(\tau)}\cdot\phi\Big)dxdt=0,
\end{aligned}
\end{equation}
for piecewise constant functions $\phi:(0,T)\to H^m(\Omega)$. 

Lemma 5, \cite{2} only requires that $h^{\prime\prime}(u)A(u)$ is positive semi-definite. The positive semi-definiteness of $h^{\prime\prime}(u)A(u)$ is included in (H2). 

\textbf{Step 2.} Relying on (\ref{inductionhyp}) and the inequality
\begin{equation}\nonumber
\begin{aligned}
h(u(w^k))-h(u(w^{k-1}))\leq h^{\prime}(u(w^k))\cdot(u(w^k)-u(w^{k-1})),
\end{aligned}
\end{equation}
we deduce that $w^k$ satisfies the discrete entropy inequality
\begin{equation}\label{entropyinequality}
\begin{aligned}
&\int_{\Omega}h(u(w^k))dx+\tau\int_{\Omega}\nabla w^k:B(w^k)\nabla w^kdx\\&+\varepsilon\tau\int_{\Omega}\Big(\sum_{|\alpha|=m}|D^{\alpha}w^k|^2+|w^k|^2\Big)dx\leq\int_{\Omega}h(u(w^{k-1}))dx.
\end{aligned}
\end{equation}

By the generalized Poincar\'e inequality, we deduce that
\begin{equation}\nonumber
\begin{aligned}
\int_{\Omega}\Big(\sum_{|\alpha|=m}|D^{\alpha}w^k|^2+|w^k|^2\Big)dx\geq C_p\|w^k\|^2_{H^m(\Omega)},
\end{aligned}
\end{equation}
where $C_p>0$ is the Poincar\'e constant. By Lemma \ref{positivesemidefinite}, we have
\begin{equation}\nonumber
\begin{aligned}
&\nabla w^k:B(w^k)\nabla w^k=\nabla u^k:H(u^k)A(u^k)\nabla u^k\geq c(|\nabla u^{k}_{1}|^2+|\nabla u^{k}_2|^2).
\end{aligned}
\end{equation}

Combining above computation results, by (\ref{entropyinequality}), we conclude that
\begin{equation}\label{konetojw}
\begin{aligned}
&\int_{\Omega}h(u^k)dx+c\tau\int_{\Omega}|\nabla u^k|^2dx+\varepsilon C_p\tau\|w^k\|^2_{H^m(\Omega)}\leq\int_{\Omega}h(u^{k-1})dx.
\end{aligned}
\end{equation}

Summing these inequalities over $k=1,...,j$ in (\ref{konetojw}), we can show that
\begin{equation}\nonumber
\begin{aligned}
&\int_{\Omega}h(u^j)dx+c\tau\sum_{k=1}^{j}\int_{\Omega}|\nabla u^k|^2dx+\varepsilon C_p\tau\sum_{k=1}^{j}\|w^k\|^2_{H^m(\Omega)}\leq\int_{\Omega}h(u^0)dx,
\end{aligned}
\end{equation}
which follows that
\begin{equation}\label{firstaubin}
\begin{aligned}
&\|u^{(\tau)}\|_{L^2(0,T;H^1(\Omega))}+\varepsilon^{1/2}\|w^{(\tau)}\|_{L^2(0,T;H^m(\Omega))}\leq C,
\end{aligned}
\end{equation}
with $C>0$ a generic constant independent of $\tau,\varepsilon$.

\textbf{Step 3.} Let us give a uniform estimate for the discrete time derivative of $u^{(\tau)}$. $A(u^{(\tau)})$ satisfies assumption (H3), then by (\ref{firstaubin}) and H\"older's inequality,
\begin{equation}\label{secondaubin}
\begin{aligned}
&\frac{1}{\tau}\Big|\int_{0}^{T}\int_{\Omega}(u^{(\tau)}-\sigma_{\tau}u^{(\tau)})\cdot\phi dxdt\Big|\leq\varepsilon\|w^{(\tau)}\|_{L^2(0,T;H^m(\Omega))}\cdot\|\phi\|_{L^2(0,T;H^m(\Omega))}\\&+\|A(u^{(\tau)})\|_{L^{\infty}(0,T;L^{\infty}(\Omega))}\cdot\|\nabla u^{(\tau)}\|_{L^2(0,T:L^2(\Omega))}\cdot\|\nabla\phi\|_{L^2(0,T;L^2(\Omega))}\leq C.
\end{aligned}
\end{equation}

By (\ref{firstaubin})-(\ref{secondaubin}), we apply the Aubin-Lions Lemma. There exists a subsequence which we do not relabel, such that as $(\tau,\varepsilon)\to 0$,
\begin{equation}\label{firstconvergence}
\begin{aligned}
u^{(\tau)}&\to u\text{ strongly in }L^2(0,T;L^2(\Omega)),\text{ and a.e. in }\Omega\times(0,T).
\\u^{(\tau)}&\rightharpoonup u\text{ weakly in }L^2(0,T;H^1(\Omega)).
\\\tau^{-1}(u^{(\tau)}-\sigma_{\tau}u^{(\tau)})&\rightharpoonup\partial_tu\text{ weakly in }L^2(0,T;H^m(\Omega)^{\prime}).\\\varepsilon w^{(\tau)}&\to 0\text{ strongly in }L^2(0,T;H^m(\Omega)).
\end{aligned}
\end{equation}

Let us denote $J=(J_1,J_2)$, with
\begin{equation}\nonumber
\begin{aligned}
&J_1=\alpha_1\nabla(u^{(\tau)}_1-u_1)+(\alpha_1-1)\big((u^{(\tau)}_1\nabla u^{(\tau)}_2-u^{(\tau)}_2\nabla u^{(\tau)}_1)-(u_1\nabla u_2-u_2\nabla u_1)\big),
\\&J_2=\alpha_2\nabla(u^{(\tau)}_2-u_2)-(\alpha_2-1)\big((u^{(\tau)}_1\nabla u^{(\tau)}_2-u^{(\tau)}_2\nabla u^{(\tau)}_1)-(u_1\nabla u_2-u_2\nabla u_1)\big).
\end{aligned}
\end{equation}

Since
\begin{equation}\label{u1u2u3difference}
\begin{aligned}
u^{(\tau)}_i\nabla u^{(\tau)}_j-u_i\nabla u_j=u^{(\tau)}_i\nabla(u^{(\tau)}_j-u_j)+(u^{(\tau)}_i-u_i)\nabla u_j,\quad 1\leq i,j\leq 2,
\end{aligned}
\end{equation}
then by (\ref{firstconvergence})-(\ref{u1u2u3difference}) and dominated convergence theorem, we can show that for every $\phi\in L^2(0,T;H^m(\Omega))$, as $(\tau,\varepsilon)\to 0$,
\begin{equation}\label{j1j3}
\begin{aligned}
\int_{0}^{T}\int_{\Omega}J\cdot\nabla\phi dxdt\to 0.
\end{aligned}
\end{equation}

(\ref{j1j3}) implies that for every $\phi\in L^2(0,T;H^m(\Omega))$, as $(\tau,\varepsilon)\to 0$, 
\begin{equation}\label{secondconvergence}
\begin{aligned}
&\int_{0}^{T}\int_{\Omega}\big(A(u^{(\tau)})\nabla u^{(\tau)}-A(u)\nabla u\big)\cdot\nabla\phi dxdt\to 0.
\end{aligned}
\end{equation}

By (\ref{firstconvergence}) and (\ref{secondconvergence}), we deduce that for every $T>0$, $\phi\in L^2(0,T;H^m(\Omega))$, we have
\begin{equation}\label{lambdalimit}
\begin{aligned}
\int_{0}^{T}\int_{\Omega}\partial_tu\cdot\phi dxdt+\int_{0}^{T}\int_{\Omega}\nabla\phi:A(u)\nabla udxdt=0.
\end{aligned}
\end{equation}

By (\ref{secondaubin}), $\partial_tu\in L^2(0,T;H^1(\Omega)^{\prime})$, and consequently, the weak formulation (\ref{lambdalimit}) also holds for every $\phi\in L^2(0,T;H^1(\Omega))$. Finally, $u\in H^1(0,T;H^1(\Omega)^{\prime})\hookrightarrow C^0([0,T];H^1(\Omega)^{\prime})$ shows that the initial condition is satisfied in $H^1(\Omega)^{\prime}$. $u$ satisfies (\ref{equation})-(\ref{initial}) in a weak sense, and we finish the proof.

\end{proof}

\section{Proof of Theorem \ref{mainresultunique}}

In the following uniqueness proof, we state an important lemma, which is Lemma 10 of \cite{4}. This lemma also has a detailed explanation in \cite{19}, so we omit the proof.
\begin{lem}\label{j2require}
If $\mu$ is an absolutely continuous measure with respect to the Lebesgue measure, and $f,g:\Omega\to[0,\infty)$ are measurable, bounded, positive functions, such that $\sqrt{f},\sqrt{g}\in H^1(\Omega;d\mu)$. Then we have
\begin{equation}\nonumber
\begin{aligned}
\int_{\Omega}|\nabla\sqrt{f+g}|^2d\mu\leq\int_{\Omega}|\nabla\sqrt{f}|^2d\mu+\int_{\Omega}|\nabla\sqrt{g}|^2d\mu.
\end{aligned}
\end{equation}
\end{lem}

Let us give the proof of Theorem \ref{mainresultunique}.

\begin{proof}

\textbf{Step 1.} Let $\zeta(t)\in H^1(\Omega)$ be the unique solution to
\begin{equation}\nonumber
\begin{aligned}
\Delta\zeta(t)=u_3(t)-v_3(t)\text{ in }\Omega,\quad\nabla\zeta\cdot\nu=0\text{ on }\partial\Omega,\quad t>0.
\end{aligned}
\end{equation}

Since $u_3-v_3\in L^2(0,T;L^2(\Omega))$, then $\zeta\in L^2(0,T;H^1(\Omega))$. As $\partial_t(u_3-v_3)\in L^2(0,T;H^1(\Omega)^{\prime})$, we have $\Delta\partial_t\zeta\in L^2(0,T;H^1(\Omega)^{\prime})$. Then we can show that for a.e. $t>0$,
\begin{equation}\nonumber
\begin{aligned}
&\frac{1}{2}\frac{d}{dt}\int_{\Omega}|\nabla\zeta|^2dx=-\langle\Delta\partial_t\zeta,\zeta\rangle=-\langle\partial_t(u_3-v_3),\zeta\rangle\\&=\alpha\int_{\Omega}\nabla(u_3-v_3)\cdot\nabla\zeta dx=-\alpha\int_{\Omega}(u_3-v_3)^2dx\leq 0.
\end{aligned}
\end{equation}

Here, $\langle\cdot,\cdot\rangle$ denotes the duality pairing of $H^1(\Omega)^{\prime}$ and $H^1(\Omega)$. It follows that
\begin{equation}\nonumber
\begin{aligned}
\int_{\Omega}|\nabla\zeta(t)|^2dx\leq\int_{\Omega}|\nabla\zeta(0)|^2dx,\quad t>0.
\end{aligned}
\end{equation}

At time $t=0$, $\Delta\zeta(0)=u_3(0)-v_3(0)=0$ in $\Omega$, thus $\nabla\zeta(0)=0$. We deduce that $|\nabla\zeta(t)|=0$ a.e. in $\Omega$, and $u_3(t)-v_3(t)=\Delta\zeta(t)=0$ in $\Omega$. In the following proof, we denote $q=u_1+u_2=v_1+v_2$.

\textbf{Step 2.} Let us denote
\begin{equation}\nonumber
\begin{aligned}
I_1&=-\int_{\Omega}(\alpha\nabla u_1+(\alpha-1)(u_1\nabla u_2-u_2\nabla u_1))\cdot\frac{\nabla u_1}{u_1+\varepsilon}dx\\&-\int_{\Omega}(\alpha\nabla u_2-(\alpha-1)(u_1\nabla u_2-u_2\nabla u_1))\cdot\frac{\nabla u_2}{u_2+\varepsilon}dx,
\end{aligned}
\end{equation}
and
\begin{equation}\nonumber
\begin{aligned}
I_2&=-\int_{\Omega}(\alpha\nabla v_1+(\alpha-1)(v_1\nabla v_2-v_2\nabla v_1))\cdot\frac{\nabla v_1}{v_1+\varepsilon}dx\\&-\int_{\Omega}(\alpha\nabla v_2-(\alpha-1)(v_1\nabla v_2-v_2\nabla v_1))\cdot\frac{\nabla v_2}{v_2+\varepsilon}dx,
\end{aligned}
\end{equation}
and
\begin{equation}\nonumber
\begin{aligned}
I_3&=\int_{\Omega}(\alpha\nabla(u_1+v_1)+(\alpha-1)(u_1\nabla u_2-u_2\nabla u_1+v_1\nabla v_2-v_2\nabla v_1))\cdot\frac{\nabla(u_1+v_1)}{u_1+v_1+2\varepsilon}dx\\&+\int_{\Omega}(\alpha\nabla(u_2+v_2)-(\alpha-1)(u_1\nabla u_2-u_2\nabla u_1+v_1\nabla v_2-v_2\nabla v_1))
\cdot\frac{\nabla(u_2+v_2)}{u_2+v_2+2\varepsilon}dx,
\end{aligned}
\end{equation}
then we have
\begin{equation}\nonumber
\begin{aligned}
\frac{d}{dt}d_{\varepsilon}(u,v)&=\sum_{i=1,2}\Big(\langle\partial_tu_i,\log(u_i+\varepsilon)\rangle+\langle\partial_tv_i,\log(v_i+\varepsilon)\rangle-\langle\partial_t(u_i+v_i),\log(\frac{u_i+v_i}{2}+\varepsilon)\rangle\Big)\\&=I_1+I_2+I_3.
\end{aligned}
\end{equation}

Since
\begin{equation}\nonumber
\begin{aligned}
u_1\nabla u_2-u_2\nabla u_1&=u_1\nabla(u_1+u_2)-(u_1+u_2)\nabla u_1=u_1\nabla q-q\nabla u_1\\&=(u_1+u_2)\nabla u_2-u_2\nabla(u_1+u_2)=q\nabla u_2-u_2\nabla q,
\end{aligned}
\end{equation}
and
\begin{equation}\nonumber
\begin{aligned}
v_1\nabla v_2-v_2\nabla v_1&=v_1\nabla(v_1+v_2)-(v_1+v_2)\nabla v_1=v_1\nabla q-q\nabla v_1\\&=(v_1+v_2)\nabla v_2-v_2\nabla(v_1+v_2)=q\nabla v_2-v_2\nabla q,
\end{aligned}
\end{equation}
which follows
\begin{equation}\nonumber
\begin{aligned}
I_1&=-\int_{\Omega}\big(\frac{|\nabla u_1|^2}{u_1+\varepsilon}+\frac{|\nabla u_2|^2}{u_2+\varepsilon}\big)\cdot(\alpha(1-q)+q)dx-(\alpha-1)\int_{\Omega}\big(\frac{u_1\nabla u_1}{u_1+\varepsilon}+\frac{u_2\nabla u_2}{u_2+\varepsilon}\big)\nabla qdx,
\end{aligned}
\end{equation}
and
\begin{equation}\nonumber
\begin{aligned}
I_2&=-\int_{\Omega}\big(\frac{|\nabla v_1|^2}{v_1+\varepsilon}+\frac{|\nabla v_2|^2}{v_2+\varepsilon}\big)\cdot(\alpha(1-q)+q)dx-(\alpha-1)\int_{\Omega}\big(\frac{v_1\nabla v_1}{v_1+\varepsilon}+\frac{v_2\nabla v_2}{v_2+\varepsilon}\big)\nabla qdx,
\end{aligned}
\end{equation}
and
\begin{equation}\nonumber
\begin{aligned}
I_3&=\int_{\Omega}\big(\frac{|\nabla(u_1+v_1)|^2}{u_1+v_1+2\varepsilon}+\frac{|\nabla(u_2+v_2)|^2}{u_2+v_2+2\varepsilon}\big)\cdot(\alpha(1-q)+q)dx\\&+(\alpha-1)\int_{\Omega}\Big(\frac{(u_1+v_1)\nabla(u_1+v_1)}{u_1+v_1+2\varepsilon}+\frac{(u_2+v_2)\nabla(u_2+v_2)}{u_2+v_2+2\varepsilon}\Big)\nabla qdx.
\end{aligned}
\end{equation}

Let us further denote
\begin{equation}\label{j1form}
\begin{aligned}
H_1&=-(\alpha-1)\int_{0}^{t}\int_{\Omega}\Big(\big(\frac{u_1\nabla u_1}{u_1+\varepsilon}+\frac{u_2\nabla u_2}{u_2+\varepsilon}\big)+\big(\frac{v_1\nabla v_1}{v_1+\varepsilon}+\frac{v_2\nabla v_2}{v_2+\varepsilon}\big)\\&-\Big(\frac{(u_1+v_1)\nabla(u_1+v_1)}{u_1+v_1+2\varepsilon}+\frac{(u_2+v_2)\nabla(u_2+v_2)}{u_2+v_2+2\varepsilon}\Big)\Big)\nabla qdxds,
\end{aligned}
\end{equation}
and
\begin{equation}\label{j2form}
\begin{aligned}
H_2&=\sum_{i=1,2}\int_{0}^{t}\int_{\Omega}\Big(\frac{|\nabla(u_i+v_i)|^2}{u_i+v_i+2\varepsilon}-\frac{|\nabla u_i|^2}{u_i+\varepsilon}-\frac{|\nabla v_i|^2}{v_i+\varepsilon}\Big)\cdot(\alpha(1-q)+q)dxds\\&=4\sum_{i=1,2}\int_{0}^{t}\int_{\Omega}\big(|\nabla\sqrt{u_i+v_i+2\varepsilon}|^2-|\nabla\sqrt{u_i+\varepsilon}|^2-|\nabla\sqrt{v_i+\varepsilon}|^2\big)\cdot(\alpha(1-q)+q)dxds\leq 0,
\end{aligned}
\end{equation}
then we can show that
\begin{equation}\label{HII}
\begin{aligned}
\int_{0}^{t}(I_1+I_2+I_3)ds=H_1+H_2.
\end{aligned}
\end{equation}

In (\ref{j2form}), we rely on Lemma \ref{j2require} to deduce that $H_2\leq 0$. For $i=1,2$, we have
\begin{equation}\nonumber
\begin{aligned}
\nabla u_i,\nabla v_i,\nabla q\in L^2(0,T;L^2(\Omega)),\quad\big|\frac{u_i}{u_i+\varepsilon}\big|\leq 1,\quad\big|\frac{v_i}{v_i+\varepsilon}\big|\leq 1,\quad\big|\frac{u_i+v_i}{u_i+v_i+2\varepsilon}\big|\leq 1,
\end{aligned}
\end{equation}
then by the dominated convergence theorem, we derive that as $\varepsilon\to 0$,
\begin{equation}\label{j1form1}
\begin{aligned}
&H_1\to 0.
\end{aligned}
\end{equation}

\textbf{Step 3.} By (\ref{j2form})-(\ref{j1form1}), since $d_{\varepsilon}(u,v)\geq 0$, $d_{\varepsilon}(u^0,v^0)=0$, we conclude that as $\varepsilon\to 0$, $d_{\varepsilon}(u,v)\to 0$. For $i=1,2$ and a.e. $(x,t)\in\Omega\times(0,T)$,
\begin{equation}\label{taylorone}
\begin{aligned}
\xi_{\varepsilon}(u_i)+\xi_{\varepsilon}(v_i)-2\xi_{\varepsilon}(\frac{u_i+v_i}{2})\to 0,\quad\text{as }\varepsilon\to 0.
\end{aligned}
\end{equation}

By the Taylor's formula, there are functions $\theta_{\varepsilon},\eta_{\varepsilon}\in[0,1]$, such that
\begin{equation}\nonumber
\begin{aligned}
&\xi_{\varepsilon}(u_i)=\xi_{\varepsilon}(\frac{u_i+v_i}{2})+\xi^{\prime}_{\varepsilon}(\frac{u_i+v_i}{2})\cdot\frac{u_i-v_i}{2}+\frac{1}{2}\xi^{\prime\prime}_{\varepsilon}\big(\theta_{\varepsilon}\cdot\frac{u_i+v_i}{2}+(1-\theta_{\varepsilon})u_i\big)\cdot\big(\frac{u_i-v_i}{2}\big)^2,\\&\xi_{\varepsilon}(v_i)=\xi_{\varepsilon}(\frac{u_i+v_i}{2})-\xi^{\prime}_{\varepsilon}(\frac{u_i+v_i}{2})\cdot\frac{u_i-v_i}{2}+\frac{1}{2}\xi^{\prime\prime}_{\varepsilon}\big(\eta_{\varepsilon}\cdot\frac{u_i+v_i}{2}+(1-\eta_{\varepsilon})v_i\big)\cdot\big(\frac{u_i-v_i}{2}\big)^2.
\end{aligned}
\end{equation}

Since $\xi^{\prime\prime}_{\varepsilon}(s)=1/(s+\varepsilon)\geq 1/2$ for $0\leq s\leq 1$, $0<\varepsilon<1$, we derive that
\begin{equation}\label{taylortwo}
\begin{aligned}
\xi_{\varepsilon}(u_i)+\xi_{\varepsilon}(v_i)-2\xi_{\varepsilon}(\frac{u_i+v_i}{2})\geq\frac{(u_i-v_i)^2}{8}.
\end{aligned}
\end{equation}

Combining (\ref{taylorone})-(\ref{taylortwo}), we can show that $u=v$, a.e. $(x,t)\in\Omega\times(0,T)$.

\end{proof}

\section{Appendix}

\subsection{Lemma 5 of \cite{2}} For integrity purpose, we give a brief proof to Lemma 5 of \cite{2}. Let $y\in L^{\infty}(\Omega;\mathbb{R}^n)$ and $\delta\in[0,1]$, and define
\begin{equation}\nonumber
\begin{aligned}
&a(w,\phi)=\int_{\Omega}\nabla\phi:B(y)\nabla wdx+\varepsilon\int_{\Omega}\Big(\sum_{|\alpha|=m}D^{\alpha}w\cdot D^{\alpha}\phi+w\cdot\phi\Big)dx,\\&F(\phi)=-\frac{\delta}{\tau}\int_{\Omega}(u(y)-u(w^{k-1}))\cdot\phi dx.
\end{aligned}
\end{equation}

Let us consider the equation
\begin{equation}\label{laxmilgram}
\begin{aligned}
&a(w,\phi)=F(\phi),\text{ for every }\phi\in H^m(\Omega;\mathbb{R}^n).
\end{aligned}
\end{equation}

The forms $a$ and $F$ are bounded on $H^m(\Omega;\mathbb{R}^n)$. The matrix $B(y)=A(u(y))(h^{\prime\prime}(u(y)))^{-1}$ is positive semi-definite since
\begin{equation}\nonumber
\begin{aligned}
z^{\mathrm{T}}B(y)z=((h^{\prime\prime}(u))^{-1}z)^{\mathrm{T}}h^{\prime\prime}(u)A(u(y))(h^{\prime\prime}(u))^{-1}z\geq 0,\text{ for every }z\in\mathbb{R}^n.
\end{aligned}
\end{equation}

For the bilinear form $a$, by the Poincar\'e inequality, we have
\begin{equation}\nonumber
\begin{aligned}
a(w,w)\geq\varepsilon\int_{\Omega}\Big(\sum_{|\alpha|=m}|D^{\alpha}w|^2+|w|^2\Big)dx\geq\varepsilon C_p\|w\|^2_{H^m(\Omega)}.
\end{aligned}
\end{equation}

We apply the Lax-Milgram lemma to obtain the existence of a unique solution $w\in H^m(\Omega;\mathbb{R}^n)\hookrightarrow L^{\infty}(\Omega;\mathbb{R}^n)$ to (\ref{laxmilgram}). This defines the fixed point operator $S:L^{\infty}(\Omega;\mathbb{R}^n)\times[0,1]\to L^{\infty}(\Omega;\mathbb{R}^n)$, $S(y,\delta)=w$, where $w$ solves (\ref{laxmilgram}).

The fixed point operator $S$ is continuous, compact and bounded. For continuity and compactness of the operator $S$, please refer to \cite{2} for a detailed explanation. We are left to prove a uniform bound for all fixed points of $S(\cdot,\delta)$ in $L^{\infty}(\Omega;\mathbb{R}^n)$. 

Let $w\in L^{\infty}(\Omega;\mathbb{R}^n)$ be such a fixed point. Then $w$ solves (\ref{laxmilgram}) with $y$ replaced by $w$. Let us choose the test function as $\phi=w$, then
\begin{equation}\nonumber
\begin{aligned}
\frac{\delta}{\tau}\int_{\Omega}(u(w)-u(w^{k-1}))\cdot wdx+\int_{\Omega}\nabla w:B(w)\nabla wdx+\varepsilon\int_{\Omega}\Big(\sum_{|\alpha|=m}|D^{\alpha}w|^2+|w|^2\Big)dx=0.
\end{aligned}
\end{equation}

The convexity of $h$ implies that $h(x)-h(y)\leq h^{\prime}(x)(x-y)$ for all $x,y\in\mathcal{D}$. Let us choose $x=u(w),y=u(w^{k-1})$ and $h^{\prime}(u(w))=w$, which follows
\begin{equation}\nonumber
\begin{aligned}
\delta\int_{\Omega}h(u(w))dx+\varepsilon\tau\int_{\Omega}\Big(\sum_{|\alpha|=m}|D^{\alpha}w|^2+|w|^2\Big)dx\leq\delta\int_{\Omega}h(u(w^{k-1}))dx.
\end{aligned}
\end{equation}

This yields an $H^m$ bound for $w$ uniform in $\delta$. The Leray-Schauder fixed point theorem shows that there exists a solution $w\in H^m(\Omega;\mathbb{R}^n)$ to (\ref{laxmilgram}), with $y$ replaced by $w$ and $\delta=1$.

\subsection{Conclusions} We have shown that the two species chemotaxis model (\ref{equation})-(\ref{initial}) has a weak solution. In \cite{20}, bounded weak solutions to a class of cross-diffusion systems were also discussed. In the proof, a class of entropy functionals were introduced to yield uniform estimates.

For some cross-diffusion systems, how to show their weak solutions exist remains a great challenge. In the global existence analysis of chemotaxis model with cross-diffusion phenomenon, we have demonstrated how to find the entropy structure of the model in this work.

For some cross-diffusion systems, their corresponding entropy structures are still hidden. Among these systems, Shigesada-Kawasaki-Teramoto type population system is the most typical one. In \cite{7}, weak solutions of the population system were shown exist, when diffusion coefficients satisfy certain assumptions. 

We still denote $u=(u_1,u_2)$, and constants $d_0,d_1,d_2>0$, $d_0\neq d_1$, $d_0\neq d_2$, with
\begin{equation}\nonumber
\begin{aligned}
\delta(u)=d_1d_2(1-u_1-u_2)+d_0(d_1u_1+d_2u_2).
\end{aligned}
\end{equation}

Three components Maxwell-Stefan equations is a system of parabolic partial differential equations in the form of (\ref{equation}), with the diffusion matrix given by
\begin{equation}\nonumber
\begin{aligned}
A(u)&=\frac{1}{\delta(u)}\left(\begin{array}{cccccc}
d_2+(d_0-d_2)u_1 & (d_0-d_1)u_1 \\
 (d_0-d_2)u_2 & d_1+(d_0-d_1)u_2 \\
\end{array}\right).
\end{aligned}
\end{equation}

For every test function $\varphi\in L^2(0,T;H^1(\Omega))$, $u=(u_1,u_2)$ satisfies that $u\in\overline{\mathcal{D}}$, and
\begin{equation}\nonumber
\begin{aligned}
&\langle\partial_tu_1,\varphi\rangle+\big\langle\frac{1}{\delta(u)}\big((d_2+(d_0-d_2)u_1)\nabla u_1+(d_0-d_1)u_1\nabla u_2\big),\nabla\varphi\big\rangle=0,
\\&\langle\partial_tu_2,\varphi\rangle+\big\langle\frac{1}{\delta(u)}\big((d_0-d_2)u_2\nabla u_1+(d_1+(d_0-d_1)u_2)\nabla u_2\big),\nabla\varphi\big\rangle=0.
\end{aligned}
\end{equation}

Once we denote
\begin{equation}\nonumber
\begin{aligned}
\beta_1=\frac{d_1}{d_0-d_1},\quad\beta_2=\frac{d_2}{d_0-d_2},
\end{aligned}
\end{equation}
we can show that
\begin{equation}\label{twomaxwell}
\begin{aligned}
\langle\partial_t(\beta_1u_1+\beta_2u_2),\varphi\rangle=\big\langle\Delta\big(\frac{\beta_1u_1}{d_1}+\frac{\beta_2u_2}{d_2}\big),\varphi\big\rangle.
\end{aligned}
\end{equation}

Let us denote $\delta_{ij}$ as Kronecker delta symbol. For one of the most typical cross-diffusion systems with degeneration phenomenon, the system is in the form of (\ref{equation}), with the diffusion matrix given by
\begin{equation}\nonumber
\begin{aligned}
A_{ij}(u)=\delta_{ij}D_iu_{n+1}+D_iu_i,\quad u_{n+1}=1-\sum_{i=1}^{n}u_i,\quad D_i>0,\quad i,j=1,...,n.
\end{aligned}
\end{equation}

For every test function $\varphi\in L^2(0,T;H^1(\Omega))$, the weak solution $u=(u_1,...,u_n)$ satisfies that $u\in\overline{\mathcal{V}}$, and
\begin{equation}\nonumber
\begin{aligned}
\langle\partial_t\big(\frac{u_i}{D_i}\big),\varphi\rangle+\langle u_{n+1}\nabla u_i-u_i\nabla u_{n+1},\nabla\varphi\rangle=0,\quad 1\leq i\leq n,
\end{aligned}
\end{equation}
which follows that
\begin{equation}\label{degeneratesum}
\begin{aligned}
\langle\partial_t\big(\sum_{i=1}^{n}\frac{u_i}{D_i}\big),\varphi\rangle=\langle\Delta\big(\sum_{i=1}^{n}u_i\big),\varphi\rangle.
\end{aligned}
\end{equation}

For semiconductor model with electron-hole scattering, the model is in the form of (\ref{equation}), with the diffusion matrix given by
\begin{equation}\nonumber
\begin{aligned}
A(u)&=\frac{1}{1+\mu_2u_1+\mu_1u_2}\left(\begin{array}{cccccc}
\mu_1(1+\mu_2u_1) & \mu_1\mu_2u_1 \\
 \mu_1\mu_2u_2 & \mu_2(1+\mu_1u_2) \\
\end{array}\right).
\end{aligned}
\end{equation}

For every test function $\varphi\in L^2(0,T;H^1(\Omega))$, the weak solution $u=(u_1,u_2)$ satisfies that
\begin{equation}\nonumber
\begin{aligned}
&\langle\partial_t\big(\frac{u_1}{\mu_1}\big),\varphi\rangle+\big\langle\frac{(1+\mu_2u_1)\nabla u_1+\mu_2u_1\nabla u_2}{1+\mu_2u_1+\mu_1u_2},\nabla\varphi\big\rangle=0,
\\&\langle\partial_t\big(\frac{u_2}{\mu_2}\big),\varphi\rangle+\big\langle\frac{\mu_1u_2\nabla u_1+(1+\mu_1u_2)\nabla u_2}{1+\mu_2u_1+\mu_1u_2},\nabla\varphi\big\rangle=0,
\end{aligned}
\end{equation}
and we can show that
\begin{equation}\label{semisum}
\begin{aligned}
\langle\partial_t\big(\frac{u_1}{\mu_1}+\frac{u_2}{\mu_2}\big),\varphi\big\rangle=\langle\Delta(u_1+u_2),\varphi\rangle.
\end{aligned}
\end{equation}

We notice that (\ref{twomaxwell})-(\ref{semisum}) are similar to each other. So far, uniqueness of weak solutions results for these cross-diffusion systems have been shown when $d_1=d_2$ in (\ref{twomaxwell}), $D_i=D_j$, $1\leq i,j\leq n$ in (\ref{degeneratesum}), and $\mu_1=\mu_2$ in (\ref{semisum}). How to show these weak solutions are unique for general coefficients cases is a major topic in the study of parabolic partial differential equations.
\\
\\

\textbf{Author contributions.} The author participated in the writing of the manuscript.
\\

\textbf{Funding.} The author acknowledges support from the National Natural Science
Foundation of China (NSFC), grant 12471206.
\\

\textbf{Data Availability.} No datasets were generated or analysed during the current study.
\\

\textbf{Competing interests} The author declares no competing interests.
\\


\end{document}